\documentclass{article}
\usepackage[latin1]{inputenc}
\usepackage{amssymb, amsmath, amsthm}
\usepackage[a4paper]{geometry}

\usepackage{footnote}
\makesavenoteenv{tabular}
\usepackage{tabularx}
\usepackage{rotating}

\def\url#1{\expandafter\string\csname #1\endcsname}

\usepackage{enumitem}
\usepackage{color}
\usepackage{amsmath}% http://ctan.org/pkg/amsmath

\theoremstyle{plain}

\newtheorem*{example*}{Example}
\newtheorem{lemma}{Lemma}[section]

\def\keywords{\vspace{.5em} % Add keywords
{\textit{Keywords}:\,\relax%
}}

\usepackage{color}
\definecolor{gr}{rgb}{0,0.5,0}
\usepackage{hyperref}

\begin{document}

\begin{center}
\Large{Eccentric pie charts and an unusual pie cutting}
\end{center}

\begin{center}
S\'andor BOZ\'OKI$^{\,\,1,2}$,
\footnotetext[1]{Laboratory on Engineering and Management
Intelligence, Research Group of Operations Research and Decision
Systems, Institute for Computer Science and Control, Hungarian
Academy of Sciences (MTA SZTAKI); Mail: 1518 Budapest, P.O.~Box 63, Hungary.}
\footnotetext[2]{Department of Operations Research and Actuarial Sciences, Corvinus University
of Budapest, Hungary
 \textit{E-mail: bozoki.sandor@sztaki.mta.hu}}
\end{center}
\begin{abstract}
The eccentric pie chart, a generalization of the traditional pie chart is introduced.
An arbitrary point is fixed within the circle and rays are drawn from it.
A sector is bounded by a pair of neighboring rays and the arc between them,
The sector's area, aimed to be equal to a given proportion,
is calculated from some well known equations in coordinate geometry.
The resulting system of polynomial and trigonometric equations can be approximated by a fully polynomial system,
once the non-polynomial functions are approximated by their Taylor series written up to the first few terms.
The roots of the polynomial system have been found by the homotopy continuation method, then
used as starting points of a Newton iteration for the original (non-polynomial) system.
The method is illustrated on a special pie cutting problem, and is applicable to a wide class
of nonlinear, non-polynomial systems.

\end{abstract}

\section{Introduction}  \label{section:intro} % \ref{section:intro}

Assume that a circular pie is to be distributed among three children such that
their shares are proportional to the children's time spent with assisting, 40, 35 and 25 minutes.
There is a three-blade pie (or pizza) cutter in the kitchen, but
it is designed for equal slices: the angle of every pair of blades is $120^{\circ}$.
Where to locate the cutter in order to get slices of area 40\%, 35\% and 25\% of the whole pie?
The answer is shown in Figure 1, the detailed solution is given in Section \ref{section:pie-cutting}.
\unitlength 1mm
\begin{center}
\begin{picture}(50,55)
\put(7,10){\resizebox{40mm}{!}{\rotatebox{0}{\includegraphics{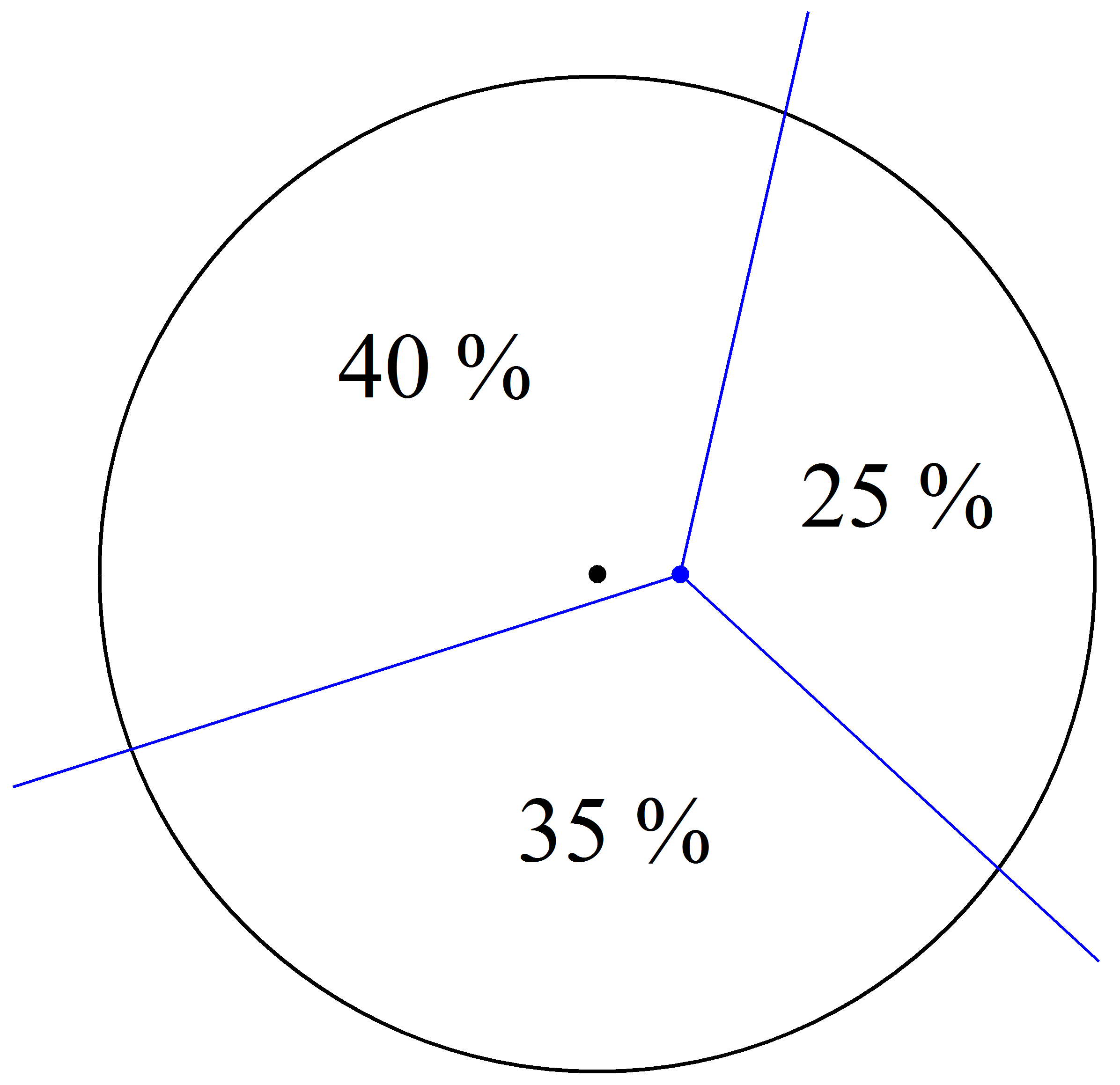}}}}
\put(-20,2){\makebox{\textbf{Figure 1.} 40-35-25\% pie-cutting with a regular 3-blade cutter }}  % \ref{figure:pie-cutting-40-35-25v2}
\end{picture}
\end{center}

Pie chart is more than two hundred years old.
According to
Spence \cite{Spence2005}
and
Tufte \cite[page 44]{Tufte2001},
Playfair \cite[Chart 2 on page 49]{Playfair1801} invented it.
The popularity of the traditional pie chart is rooted in its simplicity and efficiency in visualization.
For an arbitrary set of $n\geq2$ positive numbers $\lambda_i, i=1,\ldots,n$ such that
$\lambda_1 + \lambda_2 + \ldots + \lambda_n = 1$,
the corresponding circle sectors in the pie chart visualize the relative magnitude of numbers $\lambda_i$ in
three equivalent ways: (i) areas, (ii) central angles, (iii) arc lengths are proportional to the
numbers $\lambda_i$.

The \emph{pizza theorem} \cite{MabryDeiermann2009,Ornes2009,Upton1967} states that
the $n$-blade ($n \geq 4$ is even) equiangular cutter, wherever it is located,
halves the circle's area by summing the areas of the alternate eccentric sectors, see Figure 2.

\unitlength 1mm
\begin{center}
\begin{picture}(50,55)
\put(7,10){\resizebox{40mm}{!}{\rotatebox{0}{\includegraphics{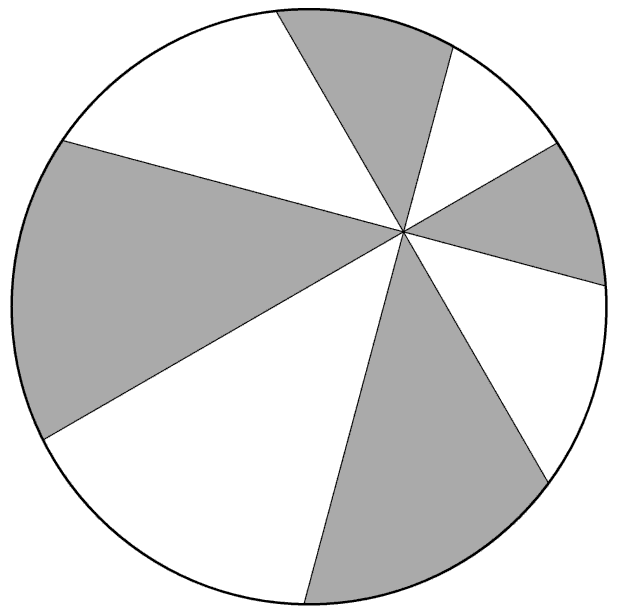}}}}
\put(-40,2){\makebox{\textbf{Figure 2.} The total area of grey
sectors is equal to the total area of white sectors }}
\put(8,-2){\makebox{(pizza theorem for $n=8$) }}
\end{picture}
\end{center}

Area-proportionality can be an important feature in case of Venn diagrams, e.g.,
for the visualization of biological lists \cite{HulsenDeVliegAlkema2008}, where
the sizes of the sets (represented by circles) as well as the overlaps correspond to the
cardinalities in the data sets.

The rest of the paper is structured as follows.
The eccentric pie chart is introduced in Section 2. It is a modification of the traditional pie chart such that
the center, the common point of the sectors is moved from the circle's center. The areas represent
the relative magnitude of numbers $\lambda_i$. The area of a sector can be calculated with the help
of a circular sector and two triangles, resulting in a system of nonlinear
polynomial and trigonometric equations. Systems of nonlinear equations are usually hard to
solve \cite[Chapter 11]{NocedalWright2006}.
If all equations are polynomial, homotopy continuation methods \cite{LeeLiTsai2008,ChenLeeLi2017} can be applied.
HOM4PS-3 \cite{ChenLeeLi2017} is used in this paper. Several geometric problems, such as
Littlewood's problem on seven mutually touching infinite
cylinders \cite{BozokiLeeRonyai2015}, or Steiner's conic problem \cite{BreidingSturmfelsTimme2019}
lead to polynomial equations. The equations of the eccentric pie chart include non-polynomials, however, an
approximation by the Taylor series, written up to the first few terms, is polynomial as
in \cite[Example 4.4]{JiWuFengLiQin2016} and \cite{Lim2004,Yalcinbas2002}.
The roots of the approximating polynomial system
can be found by the homotopy method, then they are used as starting points of a Newton iteration
for the approximated system of non-polynomial equations.

The solution of the 40-35-25\% pie-cutting problem is presented in Section 3.
The problem leads to a system of 11 polynomial equations of 11 variables. Section 4 concludes.

%Section 4 concludes with a visualization of the Pareto principle \cite{Juran1975,Pareto1906}
%(80\% of the effects come from 20\% of the causes,
% 80\% of the wealth is owned by the richest 20\% of the population, etc.) by an eccentric pie chart,
%and raises some open problems.

\newpage
\section{Eccentric pie charts}  \label{section:epc} % \ref{section:epc}

Consider the point $(x_0,y_0)$ inside the unit circle, and draw rays from it.
A sector is bounded by a pair of neighboring rays and the arc between them (Figure 3).
We focus on the area of the sectors. Unlike in case of the traditional pie chart, the angle
of the neighboring rays and the arc length are not proportional to the sector's area any more.

\unitlength 1mm
\begin{center}
\begin{picture}(50,50)
\put(0,0){\resizebox{50mm}{!}{\rotatebox{0}{\includegraphics{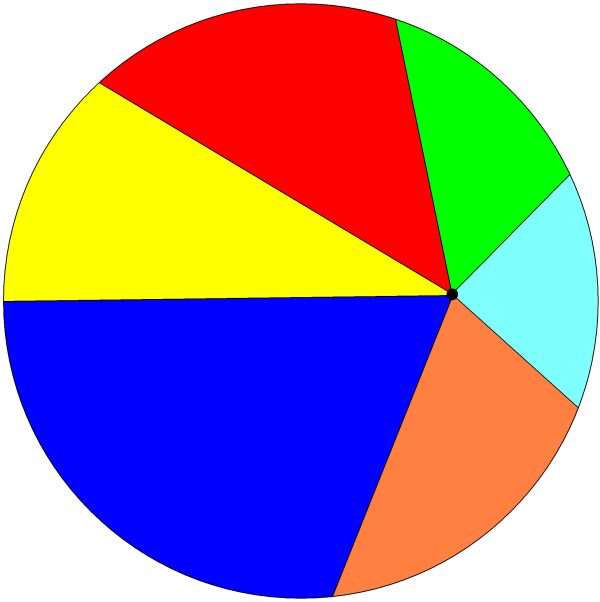}}}}
\end{picture}
\begin{center}
\textbf{Figure 3.} An eccentric pie chart
\end{center}
\end{center}

There are infinitely many eccentric pie charts representing the
same set of proportions (percentages). In fact the degree of
freedom of eccentric pie charts is 2, if rotations of the circle
are not distinguished. Once the starting point, e.g. (0,1), on the
boundary of the circle is fixed, the point $(x_0,y_0)$ can be
located anywhere inside the unit circle, there exists exactly one
eccentric pie chart representing the given set of proportions
counterclockwise, and another one clockwise. Figure 4 shows nine
eccentric pie charts ($x_0,y_0 \in \{-1/2,0,1/2\}$), all of them
visualize the 20-30-15-25-10\% counterclockwise. The first ray is
between points $(x_0,y_0)$ and (0,1). The pie chart in the middle
is the traditional pie chart.

\unitlength 1mm
\begin{center}
\begin{picture}(130,130)
\put(0,0){\resizebox{130mm}{!}{\rotatebox{0}{\includegraphics{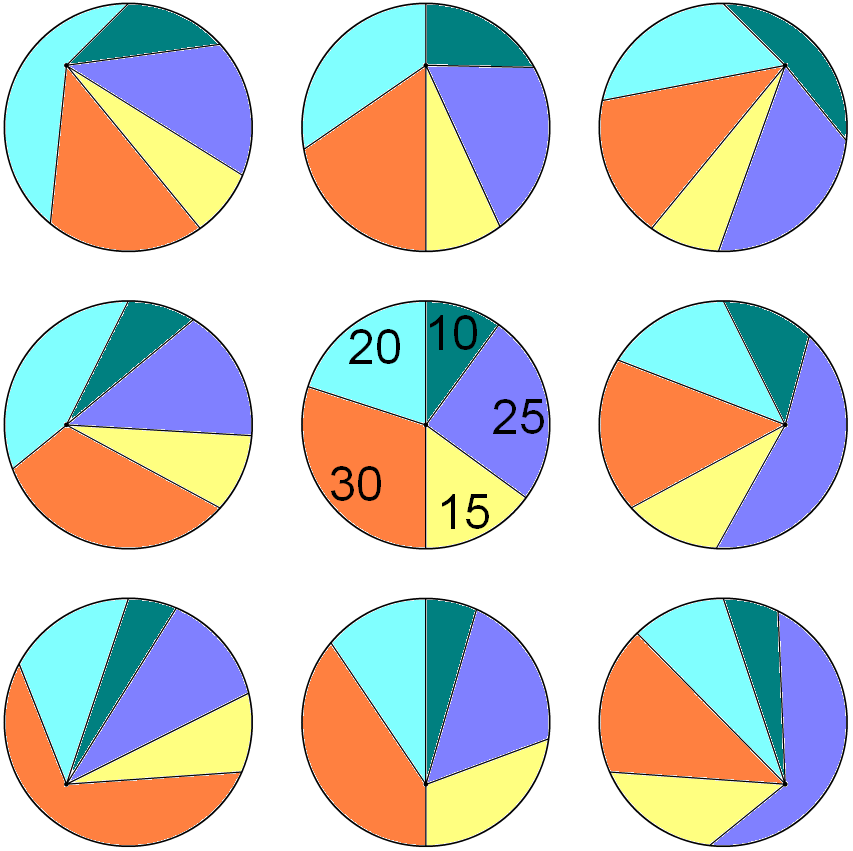}}}}
\end{picture}
\begin{center}
\textbf{Figure 4.} Eccentric pie charts with areas 20-30-15-25-10\%, where $x_0,y_0 \in \{-1/2,0,1/2\}$
\end{center}
\end{center}

The areas of eccentric sectors can be calculated from the area of a circular sector and two triangles
(Figures 5 and 6).
Let $(x_1,y_1)$ and $(x_2,y_2)$ be two points on the boundary of the unit circle centered at the origin.
The area of a circular sector is equal to $\frac{\beta}{2}$, where $\beta$ is the central angle.
It is also well known that $x_1 x_2 + y_1 y_2 = \cos \beta$.
The area of the triangle $(x_1,y_1),(x_2,y_2),(0,0)$ is $\frac{\sin \beta}{2}$, if $\beta < \pi$ (Figure 5),
and $\frac{\sin(2\pi- \beta)}{2}$, if $\beta > \pi$ (Figure 6).

\unitlength 1mm
\begin{center}
\begin{picture}(140,80)
\put(-6,2){\resizebox{66mm}{!}{\rotatebox{0}{\includegraphics{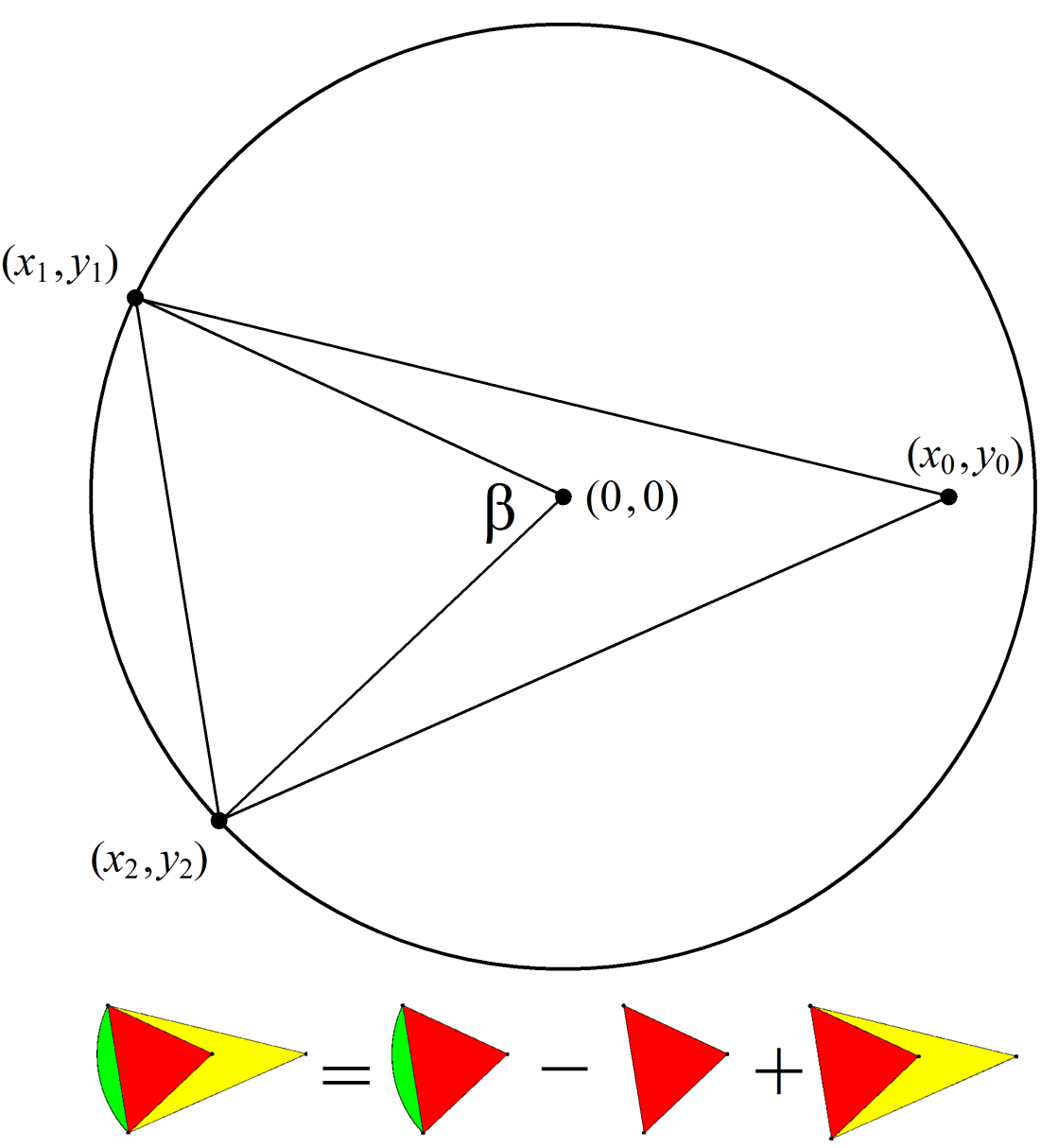}}}}
\put(80,0){\resizebox{60mm}{!}{\rotatebox{0}{\includegraphics{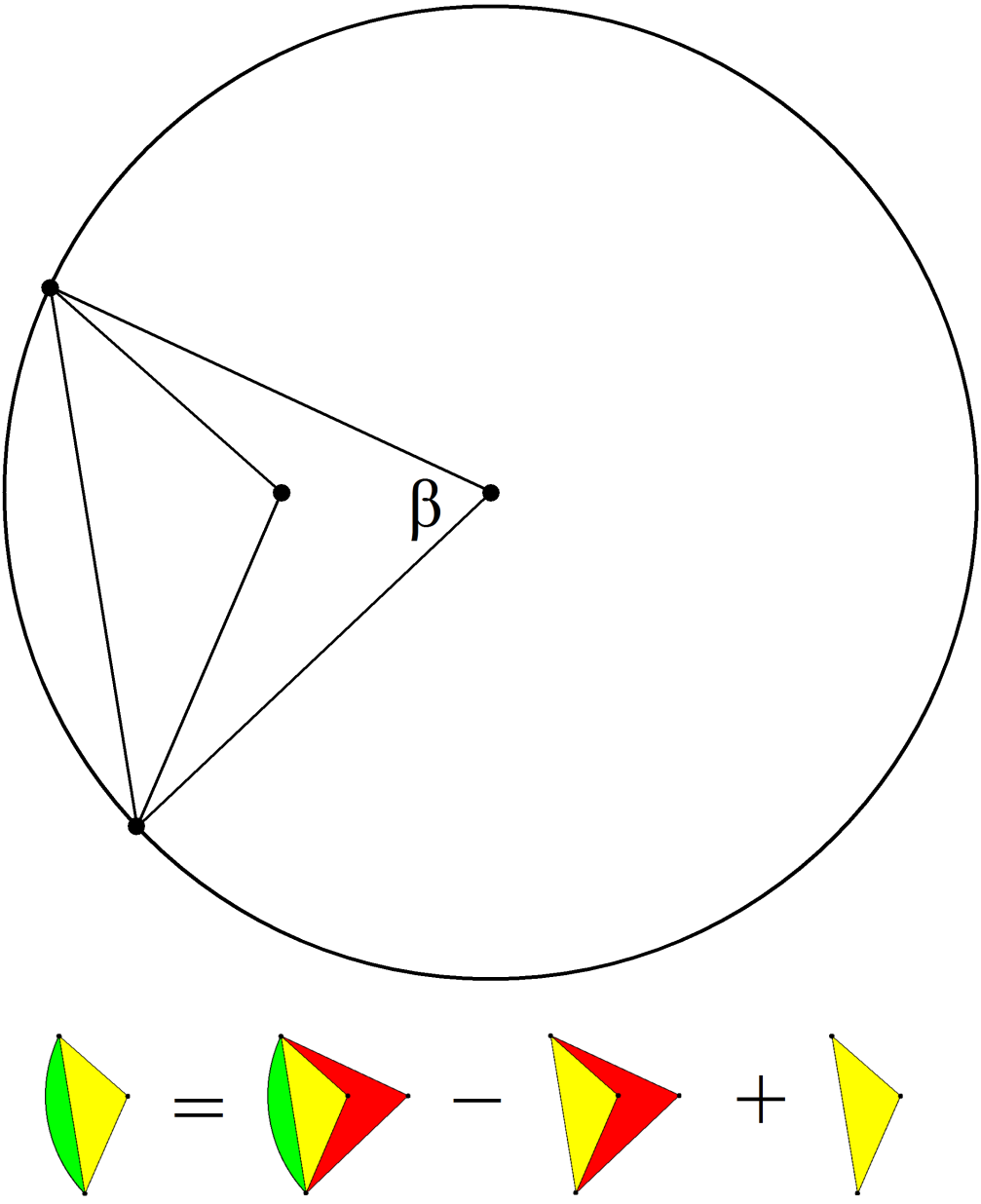}}}}
\end{picture}
\begin{center}
\textbf{Figure 5.} The calculation of the eccentric sector's area ($\beta < \pi$)
\end{center}
\end{center}

\unitlength 1mm
\begin{center}
\begin{picture}(140,80)
\put(0,0){\resizebox{60mm}{!}{\rotatebox{0}{\includegraphics{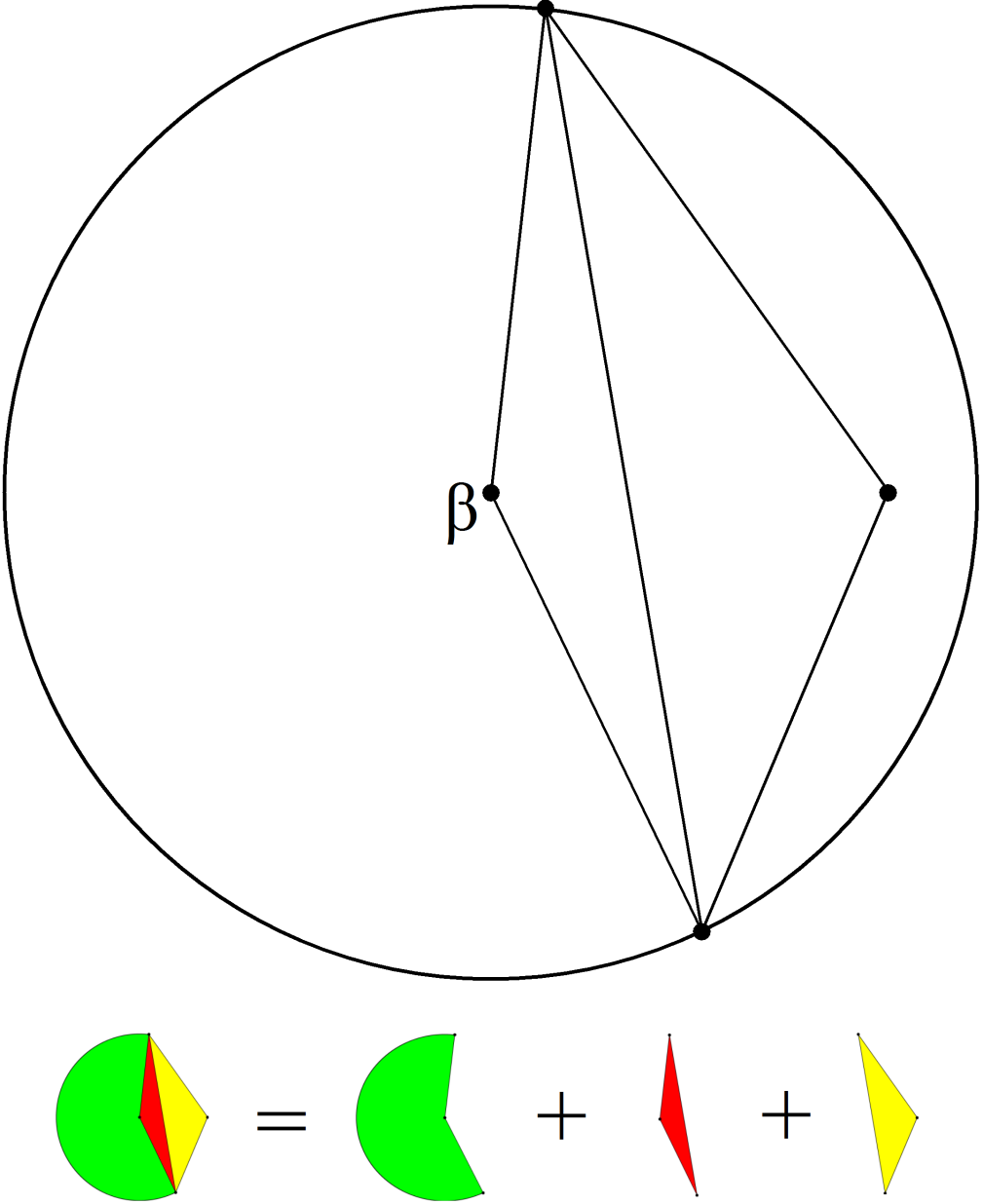}}}}
\put(80,0){\resizebox{60mm}{!}{\rotatebox{0}{\includegraphics{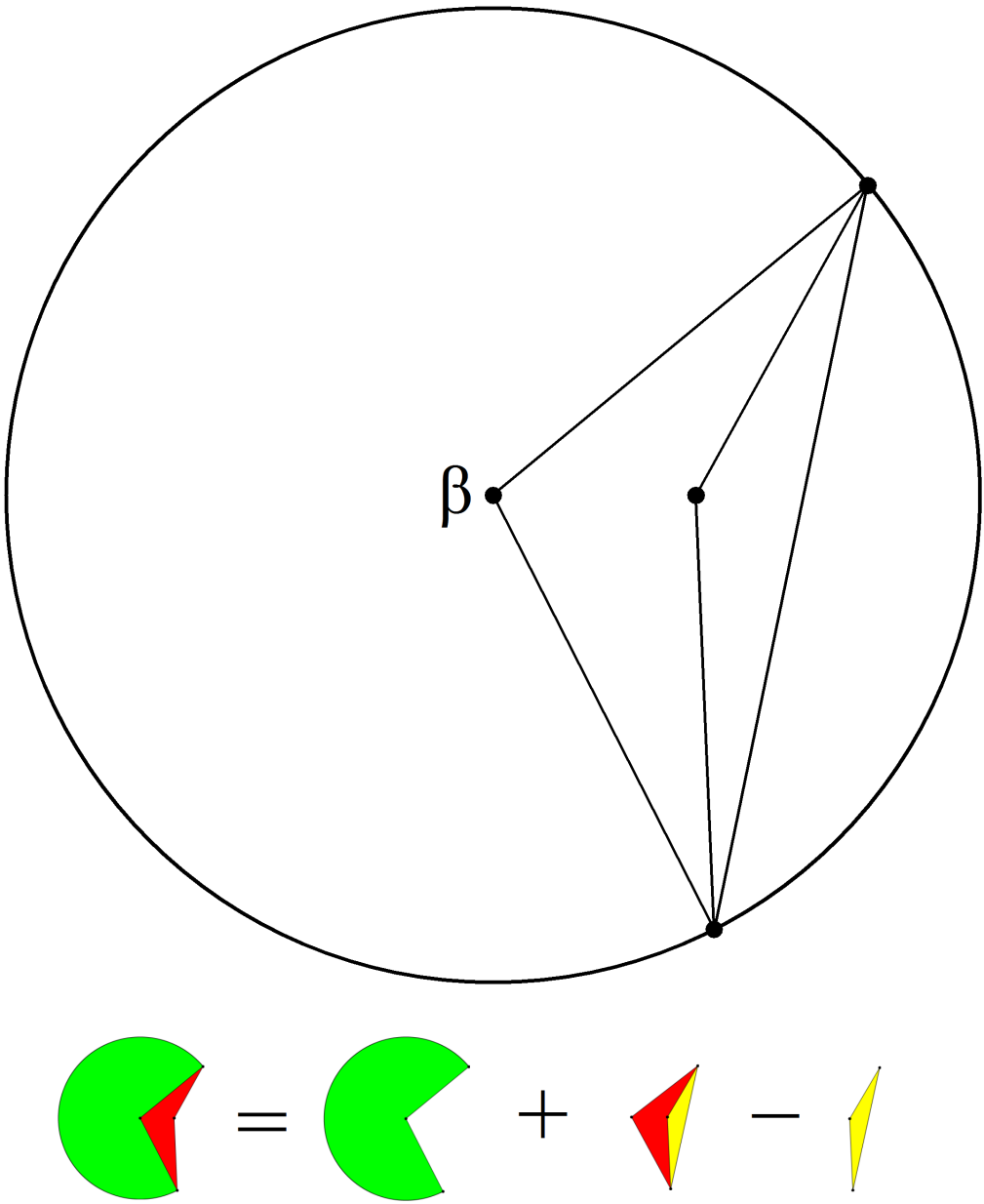}}}}
%\put(30,0){\makebox{\textbf{Figure 3.} }}
\end{picture}
\begin{center}
\textbf{Figure 6.} The calculation of the eccentric sector's area ($\beta > \pi$)
\end{center}
\end{center}

The area of a triangle can be directly calculated from the coordinates of its vertices.
\begin{lemma} (See, e.g., \cite[Problem 52 on page 34]{FineThompson1909} or
\cite[formula (4.7.2) on page 212]{Zwillinger2018} )
The area of triangle $(x_1,y_1),(x_2,y_2),(x_0,y_0)$ is equal to
\begin{equation}
\frac{1}{2}
\left|
\det\left(
\begin{matrix}
x_1 & y_1 & 1 \\
x_2 & y_2 & 1 \\
x_0 & y_0 & 1
\end{matrix}
\right)
\right|
 = \frac{1}{2}\left|
x_1 y_2 + x_2 y_0 + x_0 y_1 - x_2 y_1 - x_0 y_2 - x_1 y_0  \label{eq:TriangleArea} %(\ref{eq:TriangleArea})
\right|.
\end{equation}
\end{lemma}

Now let us prescribe that the area of the eccentric sector $(x_1,y_1),(x_2,y_2),(x_0,y_0)$
with central angle $\beta < \pi$ must be $\lambda \pi$, where $\lambda < 1$ is given.
The following equations can be written:

\begin{align}
 x_1 x_2 + y_1 y_2 = \cos \beta, \label{eq:areas-general-1} \\ %(\ref{eq:areas-general-1})
 2\left( \lambda \pi - \frac{\beta}{2} + \frac{\sin \beta}{2} \right) = | x_0(y_1-y_2) + x_1(y_2-y_0) + x_2(y_0-y_1) | , \label{eq:areas-general-2} \\ %(\ref{eq:areas-general-2})
 x_1^2 + y_1^2 = 1, \label{eq:areas-general-3} \\ %(\ref{eq:areas-general-3})
 x_2^2 + y_2^2 = 1,  \label{eq:areas-general-4}  %(\ref{eq:areas-general-4})
\end{align}

The system of equations above includes both polynomials and trigonometric expressions,
which is hard to solve. Since there exist powerful algorithms \cite{ChenLeeLi2017,LeeLiTsai2008}
for solving polynomial systems, our aim is to build a polynomial system,
which is \emph{close} to the non-polynomial system (\ref{eq:areas-general-1})-(\ref{eq:areas-general-4}).
Replace $\sin\beta$ by the new variable $s_{\beta}$, and consider the non-polynomial equation
$\beta = \arccos(x_1 x_2 + y_1 y_2) $ from (\ref{eq:areas-general-1}).

\begin{lemma} (See, e.g.,
\cite[Section 1.9.4.2 on page 50 and Section 1.9.6.5 on page 61]{Zwillinger2018}) The Taylor series of function
$\arccos$ around $a$ and $0$.
\begin{equation}
\arccos(x) = \sum\limits_{n=0}^{\infty}
\frac{{\arccos}^{(n)}(a)}{n!}(x-a)^n,  \label{eq:arccosTaylor} %(\ref{eq:arccosTaylor})
\end{equation}
especially with $a=0$
\begin{equation}
\arccos(x) = \frac{\pi}{2} - \sum\limits_{n=0}^{\infty}
\frac{(2n)!}{4^n(n!)^2(2n+1)}x^{2n+1} =
\frac{\pi}{2} - x - \frac{1}{6}x^3 - \frac{3}{40}x^5 + O(x^7), \, (|x| < 1)
\label{eq:arccosTaylor0} %(\ref{eq:arccosTaylor0})
\end{equation}
\end{lemma}

Figure 7 shows that the Taylor series around $a=0$ up to 6 terms approximates the  $\arccos$
function well if $|x| < 0.8$. If $|x|$ is close to 1, the Taylor series around $a=0.9$ up to 6 terms
provides a good approximation (Figure 8). Then one of the equations
(\ref{eq:arccosTaylor}),(\ref{eq:arccosTaylor0}) with $x = x_1 x_2 + y_1 y_2$ results in a multivariate polynomial.
The absolute values in (\ref{eq:areas-general-2}) can be eliminated by taking the squares of both sides
of the equation (false roots are possible, and they have to be filtered out).

\unitlength 1mm
\begin{center}
\begin{picture}(100,60)
\put(05,05){\resizebox{80mm}{!}{\rotatebox{0}{\includegraphics{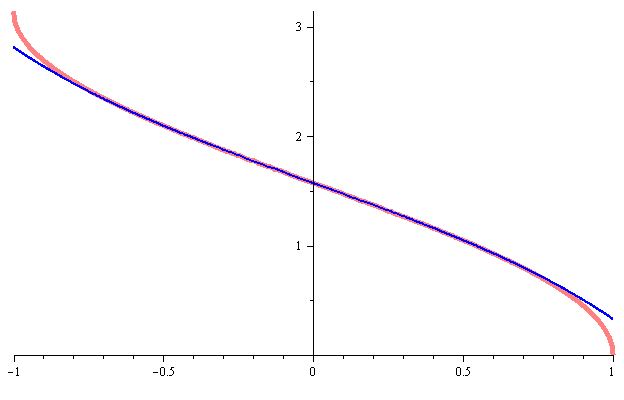}}}}
\put(81,19){\makebox{\color{blue}$\frac{\pi}{2} - x -
\frac{1}{6}x^3 - \frac{3}{40}x^5$ \color{black} }}
\put(85,12){\makebox{\color{red}$\arccos(x)$ \color{black}}}
%\put(10,5){\makebox{\textbf{Figure .} $\arccos(x)$ and $\frac{\pi}{2} - x - \frac{1}{6}x^3 - \frac{3}{40}x^5$ }}
%\put(0,5){\makebox{\textbf{Figure \label{figure:arccosTaylor0}.} $\arccos(x)$ and its Taylor series around $a=0$ up to 6 terms}} % (\ref{figure:arccosTaylor0})
%\put(20,0){\makebox{($x^2$ and $x^4$ have zero coefficients)}}
\end{picture}
\end{center}
\begin{center}
\textbf{Figure 7.} $\arccos(x)$ and its Taylor series around $a=0$ up to 6 terms\\ % (\ref{figure:arccosTaylor0})
($x^2$ and $x^4$ have zero coefficients)\\
\end{center}

\unitlength 1mm
\begin{center}
\begin{picture}(100,55)
\put(05,05){\resizebox{80mm}{!}{\rotatebox{0}{\includegraphics{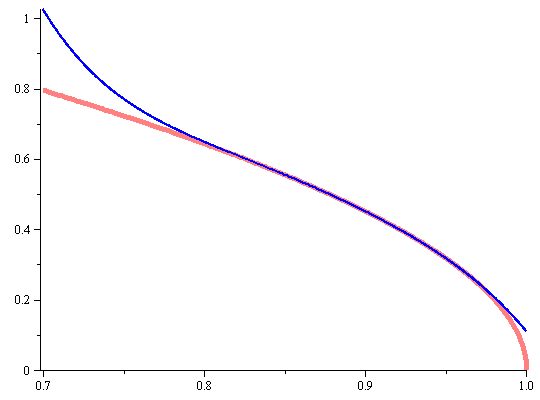}}}}
\put(21,53){\makebox{\color{blue}$0.45-2.29(x-0.9)-5.43(x-0.9)^2-27.75(x-0.9)^3$
\color{black} }}
\put(41,48){\makebox{\color{blue}$-173.84(x-0.9)^4-1218.58(x-0.9)^5$
\color{black} }} \put(85,12){\makebox{\color{red}$\arccos(x)$
\color{black}}}
%\put(0,5){\makebox{\textbf{Figure \label{figure:arccosTaylor0p9}.} $\arccos(x)$ and its Taylor series around $a=0.9$ up to 6 terms}  } % (\ref{figure:arccosTaylor0p9})
\end{picture}
\end{center}
\begin{center}
\textbf{Figure 8.} $\arccos(x)$ and its Taylor series around $a=0.9$ up to 6 terms   % (\ref{figure:arccosTaylor0p9})
\end{center}

\section{Pie cutting with a multi-blade cutter}   \label{section:pie-cutting} % \ref{section:pie-cutting}

Let us solve the 40-35-25\% pie-cutting problem from the beginning of the paper.
The geometrical problem is transformed into an algebraic one.
A system of equations is developed, then the solutions are filtered.

Let us have a unit circle with its center in the origin.
Denote by $(x_0,y_0)$ the coordinates of the 3-blade cutter's center.
We can assume that $y_0 = 0$ due to rotational symmetry.
$-1 < x_0 < 1$ is assumed in an implicit way: no equation is generated, but
after the system of equations is solved, roots not satisfying this double inequality are filtered out.

Denote by $(x_1,y_1), (x_2,y_2)$ and $(x_3,y_3)$  the coordinates of the three points, where the
boundary of the circle and the blades meet, as in Figure 9.
Let $\beta$ denote the central angle of the first eccentric sector (of area $\lambda_1 \pi = 0.4 \pi)$
bounded by line sections $[(x_0,0),(x_1,y_1)]$ and $[(x_0,0),(x_2,y_2)]$
and the arc between them.
Similarly, let $\varphi$ denote the central angle of the second eccentric sector (of area $\lambda_2 \pi = 0.35 \pi)$
bounded by line sections $[(x_0,0),(x_2,y_2)]$ and $[(x_0,0),(x_3,y_3)]$
and the arc between them. The third eccentric sector's area is $\lambda_3 \pi = 0.25 \pi$, and its central
angle is $2\pi-\beta-\varphi.$

\unitlength 1mm
\begin{center}
\begin{picture}(70,60)
\put(7,0){\resizebox{60mm}{!}{\rotatebox{0}{\includegraphics{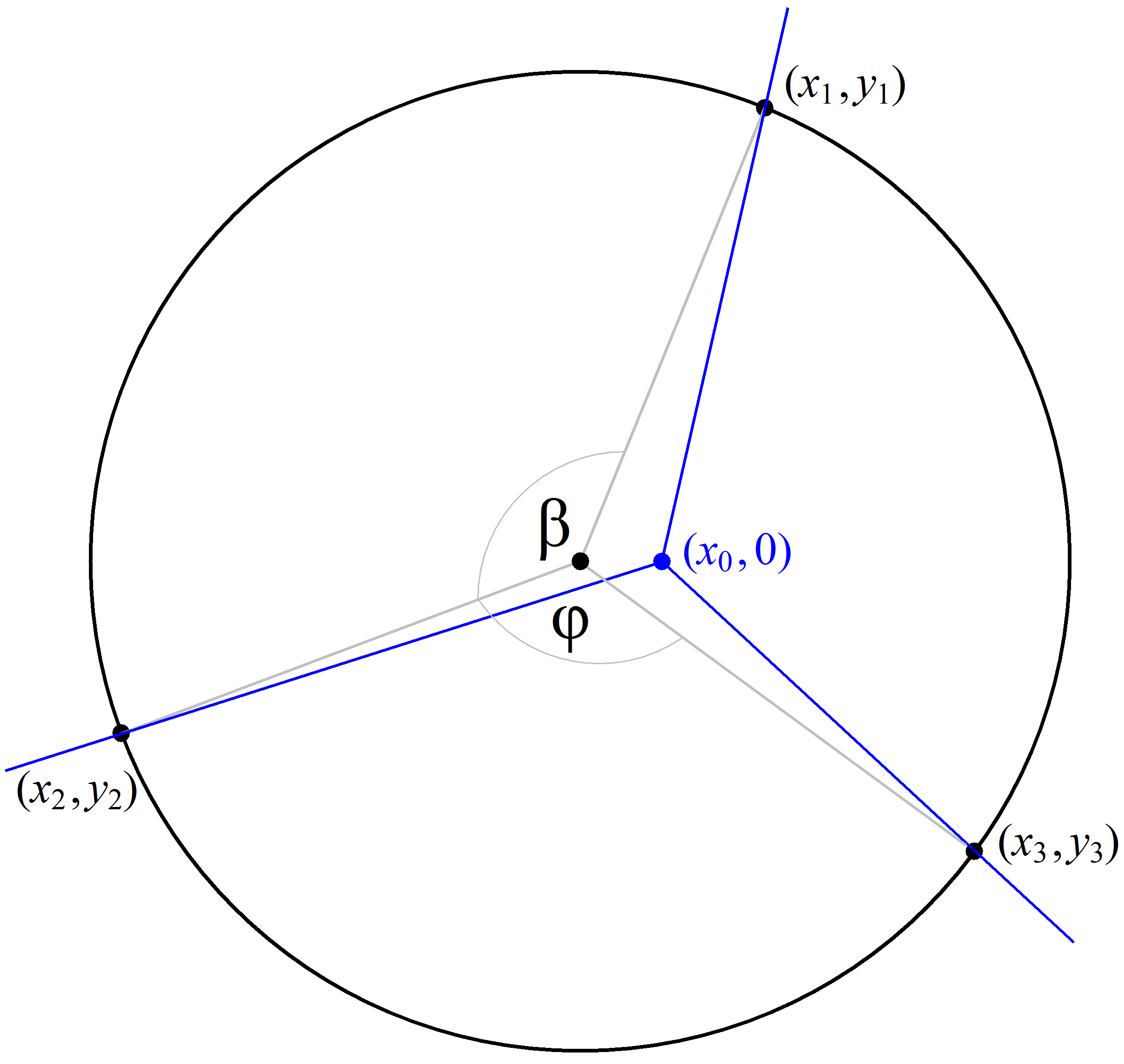}}}}
\end{picture}
\end{center}
\begin{center}
\textbf{Figure 9.} Pie-cutting with a regular 3-blade cutter
\end{center}

Following Section  \ref{section:epc} and including the regularity of 3-blade cutter (pairwise angles are $\frac{2\pi}{3}$)
we have the following equations:
\begin{align}
 (x_1 x_2 + y_1 y_2)^2 + \sin^2\beta = 1, \label{eq:pie-cutting-40-35-25-exact-01} \\ %(\ref{eq:pie-cutting-40-35-25-exact-01})
 (x_2 x_3 + y_2 y_3)^2 + \sin^2\varphi = 1, \label{eq:pie-cutting-40-35-25-exact-02} \\ %(\ref{eq:pie-cutting-40-35-25-exact-02})
   2 \lambda_1 \pi -  \beta  + \sin \beta   = | x_0(y_1-y_2) + x_1(y_2-y_0) + x_2(y_0-y_1) | , \label{eq:pie-cutting-40-35-25-exact-03} \\ %(\ref{eq:pie-cutting-40-35-25-exact-03})
   2 \lambda_2 \pi - \varphi + \sin \varphi = | x_0(y_2-y_3) + x_2(y_3-y_0) + x_3(y_0-y_2) | , \label{eq:pie-cutting-40-35-25-exact-04} \\ %(\ref{eq:pie-cutting-40-35-25-exact-04})
 x_1^2 + y_1^2 = 1, \label{eq:pie-cutting-40-35-25-exact-05} \\ %(\ref{eq:pie-cutting-40-35-25-exact-05})
 x_2^2 + y_2^2 = 1,  \label{eq:pie-cutting-40-35-25-exact-06} \\  %(\ref{eq:pie-cutting-40-35-25-exact-06})
 x_3^2 + y_3^2 = 1,  \label{eq:pie-cutting-40-35-25-exact-07} \\ %(\ref{eq:pie-cutting-40-35-25-exact-07})
\hspace*{-25mm} (x_1-x_0)(x_2-x_0) + (y_1-y_0)(y_2-y_0)  = \cos\left(\frac{2\pi}{3}\right)\sqrt{ [(x_1-x_0)^2+(y_1-y_0)^2][(x_2-x_0)^2+(y_2-y0)^2]}, \label{eq:pie-cutting-40-35-25-exact-08} \\ %(\ref{eq:pie-cutting-40-35-25-exact-08})
\hspace*{-25mm} (x_2-x_0)(x_3-x_0) + (y_2-y_0)(y_3-y_0)  = \cos\left(\frac{2\pi}{3}\right)\sqrt{ [(x_2-x_0)^2+(y_2-y_0)^2][(x_3-x_0)^2+(y_4-y0)^2]} \label{eq:pie-cutting-40-35-25-exact-09} %(\ref{eq:pie-cutting-40-35-25-exact-09})
\end{align}
where $\lambda_1 = 0.4, \lambda_2 = 0.35, y_0 = 0$ and $\cos\left(\frac{2\pi}{3}\right) = -\frac{1}{2}.$

Introduce variables $s_{\beta},s_{\varphi}$ to replace $\sin\beta$ and $\sin\varphi$, respectively.
In order to avoid absolute values in (\ref{eq:pie-cutting-40-35-25-exact-03})-(\ref{eq:pie-cutting-40-35-25-exact-04}),
take the squares of both sides. This step may bring false solutions, and we will see that it does so indeed.
The square root in (\ref{eq:pie-cutting-40-35-25-exact-08})-(\ref{eq:pie-cutting-40-35-25-exact-09})
can be eliminated likewise, with another possibility to have false solutions.
Equations (\ref{eq:pie-cutting-40-35-25-exact-01})-(\ref{eq:pie-cutting-40-35-25-exact-04})
and  (\ref{eq:pie-cutting-40-35-25-exact-08})-(\ref{eq:pie-cutting-40-35-25-exact-09})
are replaced by polynomial equations
\begin{align}
 (x_1 x_2 + y_1 y_2)^2 + s^2_{\beta} = 1, \label{eq:pie-cutting-40-35-25-poly-01} \\ %(\ref{eq:pie-cutting-40-35-25-poly-01})
 (x_2 x_3 + y_2 y_3)^2 + s^2_{\varphi} = 1, \label{eq:pie-cutting-40-35-25-poly-02} \\ %(\ref{eq:pie-cutting-40-35-25-poly-02})
 \left( 2 \lambda_1 \pi  -  \beta  + s_{\beta} \right)^2 = [x_0(y_1-y_2) + x_1(y_2-y_0) + x_2(y_0-y_1)]^2 , \label{eq:pie-cutting-40-35-25-poly-03} \\ %(\ref{eq:pie-cutting-40-35-25-poly-03})
 \left( 2 \lambda_2 \pi  - \varphi + s_{\varphi} \right)^2 = [x_0(y_2-y_3) + x_2(y_3-y_0) + x_3(y_0-y_2)]^2 , \label{eq:pie-cutting-40-35-25-poly-04} \\ %(\ref{eq:pie-cutting-40-35-25-poly-04})
\hspace*{-25mm} [(x_1-x_0)(x_2-x_0) + (y_1-y_0)(y_2-y_0)]^2  = \frac{1}{4}  [(x_1-x_0)^2+(y_1-y_0)^2][(x_2-x_0)^2+(y_2-y0)^2], \label{eq:pie-cutting-40-35-25-poly-08} \\ %(\ref{eq:pie-cutting-40-35-25-poly-08})
\hspace*{-25mm} [(x_2-x_0)(x_3-x_0) + (y_2-y_0)(y_3-y_0)]^2  = \frac{1}{4}  [(x_2-x_0)^2+(y_2-y_0)^2][(x_3-x_0)^2+(y_4-y0)^2] \label{eq:pie-cutting-40-35-25-poly-09} %(\ref{eq:pie-cutting-40-35-25-poly-09})
\end{align}
Finally approximate the equations
\begin{align}
 \cos\beta = x_1 x_2 + y_1 y_2 \label{eq:pie-cutting-40-35-25-exact-10} \\ %(\ref{eq:pie-cutting-40-35-25-exact-10})
 \cos\varphi = x_2 x_3 + y_2 y_3 \label{eq:pie-cutting-40-35-25-exact-11}  %(\ref{eq:pie-cutting-40-35-25-exact-11})
\end{align}
or, equivalently, $\beta = \arccos(x_1 x_2 + y_1 y_2)$ and $\varphi = \arccos(x_2 x_3 + y_2 y_3)$
by the Taylor series of function $\arccos$ around zero up to the fifth power as in (\ref{eq:arccosTaylor0}).
Here we expect that angles $\beta$ and $\varphi$ are between $\arccos(0.8) \approx 0.6435 \approx 37^{\circ}$
and $\pi - \arccos(0.8) \approx 2.498 \approx 143^{\circ}$. Would this assumption fail, we can try the other
Taylor series around 0.9 as in (\ref{eq:arccosTaylor}) and Figure 8.

\begin{align}
\beta  = \frac{\pi}{2} - (x_1 x_2 + y_1 y_2)^ - \frac{1}{6}(x_1 x_2 + y_1 y_2)^3 - \frac{3}{40}(x_1 x_2 + y_1 y_2)^5 , \label{eq:pie-cutting-40-35-25-sinbeta-Taylor} \\ %(\ref{eq:pie-cutting-40-35-25-sinbeta-Taylor})
\varphi = \frac{\pi}{2} - (x_2 x_3 + y_2 y_3)^ - \frac{1}{6}(x_2 x_3 + y_2 y_3)^3 - \frac{3}{40}(x_2 x_3 + y_2 y_3)^5 , \label{eq:pie-cutting-40-35-25-sinphi-Taylor}  %(\ref{eq:pie-cutting-40-35-25-sinphi-Taylor})
\end{align}
The polynomial system
(\ref{eq:pie-cutting-40-35-25-exact-05})-(\ref{eq:pie-cutting-40-35-25-exact-07}),(\ref{eq:pie-cutting-40-35-25-poly-01})-(\ref{eq:pie-cutting-40-35-25-poly-09}),(\ref{eq:pie-cutting-40-35-25-sinbeta-Taylor})-(\ref{eq:pie-cutting-40-35-25-sinphi-Taylor}) has 11 equations and
11 variables:\linebreak $x_0, x_1, y_1, x_2, y_2, x_3, y_3,
\beta,\varphi, s_{\beta}, s_{\varphi} $. Homotopy algorithm
HOM4PS-3 \cite{ChenLeeLi2017} found 28224 roots, 720 of them are
real. However, most of them are false solutions to the geometric
problem, due to several reasons. Many roots do not fulfil $ |x_0|,
|s_{\beta}| , |s_{\varphi}| \leq 1$. Some roots satisfy
(\ref{eq:pie-cutting-40-35-25-poly-03})/(\ref{eq:pie-cutting-40-35-25-poly-04})
but not
(\ref{eq:pie-cutting-40-35-25-exact-03})/(\ref{eq:pie-cutting-40-35-25-exact-03}),
or, similarly, satisfy
(\ref{eq:pie-cutting-40-35-25-poly-08})/(\ref{eq:pie-cutting-40-35-25-poly-09})
but not
(\ref{eq:pie-cutting-40-35-25-exact-08})/(\ref{eq:pie-cutting-40-35-25-exact-09}).
Furthermore, $s_{\beta} \not\approx  \sin \beta $, but instead of
that $s_{\beta} \approx -\sin\beta. $

After all 4 solutions remain to use as starting points of a Newton-iteration
for the system of equations
(\ref{eq:pie-cutting-40-35-25-exact-05})-(\ref{eq:pie-cutting-40-35-25-exact-07}),
(\ref{eq:pie-cutting-40-35-25-poly-01})-(\ref{eq:pie-cutting-40-35-25-poly-09}),
(\ref{eq:pie-cutting-40-35-25-exact-10})-(\ref{eq:pie-cutting-40-35-25-exact-11}).
Maple's \verb"fsolve" refines the solution with an arbitrary accuracy.
However we observed that the roots calculated from the polynomial system
were already within an error of 0.002 for all variables, which is due to that
the Taylor series provided a good polynomial approximation of the function $\arccos$.

The four solutions are essentially the same: one solution, given below up to 10 correct digits,
\begin{align*}
x_0 &=                        0.164641996                                  \\
x_1 &=                        0.375176778                                  \\
y_1 &=                        0.926953281                                  \\
x_2 &=                      - 0.939722783                                  \\
y_2 &=                      - 0.341937259                                  \\
x_3 &=                        0.805164109                                  \\
y_3 &=                      - 0.593052069                                  \\
\beta &=                      2.304361451 \approx 132^{\circ}              \\
\varphi &=                    2.157770813 \approx 123.6^{\circ}            \\
s_{\beta} &= \sin \beta     = 0.742792198                                  \\
s_{\varphi} &= \sin \varphi = 0.832620150
\end{align*}
has already been shown in Figures 1 and 9,
the others are its reflections on the vertical and/or horizontal axes.

Note that $\lambda_i$ ($i=1,2,3$) cannot be arbitrary in this problem. Even if the
the center $(x_0,y_0)$ is located at (1,0), i.e., on the circle's border, the largest sector's area
is at most $\pi - 2 \left( \frac{\pi}{6} - \frac{\sqrt{3}}{4} \right) \approx 0.9423 \pi.$

\section{Conclusions}  \label{section:conclusions} % \ref{section:conclusions}
The eccentric pie chart, a generalization of the traditional pie chart,
visualizes the proportions is infinitely many ways, with degree of freedom 2.
This variety suggests that additional information can also be taken into consideration, such as
the constraints on the angles of the rays in the 40-35-25\% pie-cutting problem.

The method of solving nonlinear, non-polynomial systems through a polynomial approximation, presented in Section 2
and illustrated on an example in Section 3, seems applicable in larger systems, too.
The replacement of a non-polynomial function by its
Taylor series has practical limitations, because a polynomial system
can also be hopelessly hard to solve.


\begin{thebibliography}{99}



\bibitem{BozokiLeeRonyai2015}
Boz\'{o}ki, S., Lee, T.L., R\'{o}nyai, L. (2015):
Seven mutually touching infinite cylinders,
Computational Geometry: Theory and Applications 48(2) pp.~87--93
% DOI 10.1016/j.comgeo.2014.08.007

\bibitem{BreidingSturmfelsTimme2019}
Breiding, P., Sturmfels, B., Timme, S. (2019):
3264 Conics in a Second,
https://arxiv.org/abs/1902.05518
% DOI

\bibitem{ChenLeeLi2017}
%Tianran Chen
%Tsung-Lin Lee
%Tien-Yien Li
Chen, T., Lee, T.-L., Li, T.-Y. (2017):
Mixed cell computation in Hom4PS-3
Journal of Symbolic Computation
79(3) pp.~516--534
% DOI 10.1016/j.jsc.2016.07.017
% https://www.sciencedirect.com/science/article/pii/S0747717116300542

\bibitem{FineThompson1909}
Fine, H.B., Thompson, H.D. (1909):
Coordinate Geometry,
The Macmillan Company, New York.
% Problem 52 on page 34


\bibitem{HulsenDeVliegAlkema2008}
% Tim Hulsen
% Jacob de Vlieg
% Wynand Alkema
Hulsen, T., de Vlieg, J., Alkema, W. (2008):
BioVenn -- a web application for the comparison and visualization of biological lists
using area-proportional Venn diagrams,
BMC Genomics 9 paper id.~488
% DOI 10.1186/1471-2164-9-488
% https://bmcgenomics.biomedcentral.com/articles/10.1186/1471-2164-9-488


\bibitem{JiWuFengLiQin2016}
% Zhenyi Ji
% Wenyuan Wu
% Yong Feng
% Yi Li
% Xiao Lin Qin
Ji, Z., Wu, W., Feng, Y., Li, Y., Qin, X.L. (2016):
Numerical Method for Real Root Isolation of Semi-Algebraic System and Its Applications,
Journal of Computational and Theoretical Nanoscience,
13(1) pp.~803--811
% DOI 10.1166/jctn.2016.4878
% https://www.ingentaconnect.com/content/asp/jctn/2016/00000013/00000001/art00118;jsessionid=15m3b5uq1ck3w.x-ic-live-01


%\bibitem{Juran1975}
%Juran, J.M. (1975):
%The Non-Pareto Principle; Mea Culpa.
%Quality Progress, pp.~8--9.
%%In: Juran on Quality by Design:
%%The New Steps for Planning Quality Into Goods and Services,
%%The Free Press, New York, pp.~68--71.
%Reprinted in:
%Juran, J.M., De Feo, J.A. (Eds.): Juran's Quality Handbook -- The Complete Guide to
%Performance Excellence, 6th Edition, McGraw-Hill, 2010, pp.~1021--1024
%% DOI 10.1.1.521.6224
%% https://www.projectsmart.co.uk/docs/the-non-pareto-principle.pdf
%% http://qpr.buaa.edu.cn/docs/20150407190945567926.pdf#page=1044

\bibitem{LeeLiTsai2008}
Lee, T.L., Li, T.Y., Tsai, C.H. (2008):
HOM4PS-2.0: a software package for solving polynomial systems by the polyhedral homotopy continuation method.
Computing
83, pp.~109--133
% DOI 10.1007/s00607-008-0015-6
% https://link.springer.com/article/10.1007%2Fs00607-008-0015-6


\bibitem{Lim2004}
% Teik-Cheng Lim
Lim, T.-C. (2004):
Application of Maclaurin series in relating interatomic potential functions: A review.
Journal of Mathematical Chemistry
36(2) pp.~147--160
% DOI 10.1023/B:JOMC.0000038772.74111.05
% https://link.springer.com/article/10.1023/B:JOMC.0000038772.74111.05


\bibitem{MabryDeiermann2009}
Mabry, R., Deiermann, P. (2009):
Of Cheese and Crust: A Proof of the Pizza Conjecture and Other Tasty Results.
The American Mathematical Monthly 116(5) pp.~423--438
% DOI 10.1080/00029890.2009.11920956
% https://www.jstor.org/stable/40391118
% https://www.tandfonline.com/doi/abs/10.1080/00029890.2009.11920956

\bibitem{NocedalWright2006}
Nocedal, J., Wright, S.J. (2006):
Numerical Optimization, 2nd edition, Springer
% Jorge Nocedal, Stephen J. Wright
% DOI 10.1007/978-0-387-40065-5_11     Chapter 11: Nonlinear Equations
% https://link.springer.com/book/10.1007/978-0-387-40065-5

%\bibitem{Pareto1906}
%Pareto, V. (1906)  Manuale di Economia Politica.
%Societa Editrice Libraria, Milano
%% https://archive.org/details/manualedieconomi00pareuoft


\bibitem{Ornes2009}
% Stephen Ornes
Ornes, S. (2009) : As easy as pie
New Scientist
204(2738) pp.~48--50
% DOI 10.1016/S0262-4079(09)63265-6
% https://www.sciencedirect.com/science/article/pii/S0262407909632656

\bibitem{Playfair1801}
Playfair, W. (1801):
The statistical breviary.
London: T.~Bensley


\bibitem{Spence2005}
Spence, I. (2005) No Humble Pie: The Origins and Usage of a statistical Chart.
Journal of Educational and Behavioral Statistics 30(4) pp.~353--368
% DOI 10.3102/10769986030004353
% http://www.psych.utoronto.ca/users/spence/Spence%202005.pdf
% http://journals.sagepub.com/doi/abs/10.3102/10769986030004353
% https://www.jstor.org/stable/pdf/3701294.pdf


\bibitem{Tufte2001}
Tufte, E.R. (2001):
The visual display of quantitative information.
2nd edition, Cheshire, Connecticut: Graphics Press.

\bibitem{Upton1967}
Upton, L.J. (1967): Problem 660,
Mathematics Magazine 40(3) p.~163
% https://www.jstor.org/stable/2688484

\bibitem{Yalcinbas2002}
% Salih Yalçinbas
Yal\c{c}inba\c{s}, S. (2002):
Taylor polynomial solutions of nonlinear Volterra-Fredholm integral equations.
Applied Mathematics and Computation
127(2-3) pp.~195--206
% DOI 10.1016/S0096-3003(00)00165-X
% https://www.sciencedirect.com/science/article/pii/S009630030000165X

\bibitem{Winkler2008a}
Winkler, P. (2008): Puzzled: Circular food,
Communications of the ACM 51(11) p.~112
% DOI 10.1145/1400214.1413439
% https://dl.acm.org/citation.cfm?id=1413439

\bibitem{Winkler2008b}
Winkler, P. (2008): Puzzled Solutions and Sources,
Communications of the ACM 51(12) p.~118
% DOI 10.1145/1409360.1409383
% https://dl.acm.org/citation.cfm?id=1409383

\bibitem{Zwillinger2018}
Zwillinger, D. (Editor) (2018):
CRC Standard Mathematical Tables and Formulas,
33rd edition, CRC Press,

\end{thebibliography}
\end{document}